\documentclass[english]{amsart}
\usepackage{babel}
\usepackage{amssymb}
\usepackage{amsmath}
\usepackage{enumerate}
\begin{document}
\def\R{\mathbb R}
\def\N{\mathbb N}
\def\Z{\mathbb Z}
\def\C{\mathbb C}
\def\Co{\mathcal{C}}
\def\F{\mathcal{F}}
\def\S{\mathcal{S}}
\def\1{\mathbf{1}}
\newtheorem{th-def}{Theorem-Definition}[section]
\newtheorem{theo}{Theorem}[section]
\newtheorem{lemm}[theo]{Lemma}
\newtheorem{prop}[theo]{Proposition}
\newtheorem{defi}[theo]{Definition}
\newtheorem{cor}[theo]{Corollary}
\def\dem{\noindent \textbf{Proof: }}
\def\rem{\indent \textsc{Remark. }}
\def\rems{\indent \textsc{Remarks. }}
\def\fin{  $\square$}
\def\ep{\varepsilon}
\def\uep{u^{\varepsilon}}
\def\Wep{W^{\varepsilon}}
\def\qep{a^{\varepsilon}}
\def\Sep{S^{\varepsilon}}
\def\Qep{Q^{\varepsilon}}
\def\Zep{Z^{\varepsilon}}
\def\vep{v^{\varepsilon}}
\def\wep{w^{\varepsilon}}
\def\fep{f^{\varepsilon}}
\def\gep{g^{\varepsilon}}
\def\hep{H^{\varepsilon}}
\def\aep{\alpha_{\varepsilon}}
\def\ph{\phi^{\ep}}

\numberwithin{equation}{section}
\title{High frequency analysis of Helmholtz equations: case of two point sources}
\author{Elise Fouassier}
\maketitle

{\small
\textbf{Abstract}. We derive the high frequency limit of the Helmholtz equation with source term when the source is the sum of two point sources. We study it in terms of Wigner measures (quadratic observables). We prove that the Wigner measure associated with the solution satisfies a Liouville equation with, as source term, the sum of the source terms that would be created by each of the two point sources taken separately. The first step, and main difficulty, in our study is the obtention of uniform estimates on the solution. Then, from these bounds, we derive the source term in the Liouville equation together with the radiation condition at infinity satisfied by the Wigner measure. \\

}
\\

\section{Introduction}

In this article, we are interested in the analysis of the high frequency limit of the following Helmholtz equation
\begin{equation}\label{h}
-i\frac{\aep}{\ep}\uep+\Delta\uep+\frac{n(x)^2}{\ep^2}\uep
=\Sep(x),\qquad x\in\R^3
\end{equation}
with $$\Sep(x)=\Sep_0(x)+\Sep_1(x)=\frac{1}{\ep^3}S_0\Big(\frac{x}{\ep}\Big)
+\frac{1}{\ep^3}S_1\Big(\frac{x-q_1}{\ep}\Big)$$ 
where $q_1$ is a point in $\R^3$ different from the origin.\\
In the sequel, we assume that the refraction index $n$ is constant, $n(x)\equiv 1$.\\
\indent The equation (\ref{h}) modelizes the propagation of a source wave in a medium with scaled refraction index $n(x)^2/\ep^2$. There, the small positive parameter $\ep$ is related to the frequency $\omega=\frac{1}{2\pi\ep}$ of $\uep$. In this paper, we study the high frequency limit, i.e. the asymptotics $\ep\to 0$. We assume that the regularizing parameter $\aep$ is positive, with $\aep\to 0$ as $\ep\to 0$. The positivity of $\aep$ ensures the existence and uniqueness of a solution $\uep$ to the Helmholtz equation (\ref{h}) in $L^2(\R^3)$ for any $\ep>0$.\\
\indent The source term $\Sep$ models a source signal that is the sum of two source signals concentrating respectively close to the origin and close to the point $q_1$ at the scale $\ep$. The concentration profiles $S_0$ and $S_1$ are given functions. Since $\ep$ is also the scale of the oscillations dictated by the Helmholtz operator $\Delta+\frac{1}{\ep^2}$, resonant interactions can occur between these oscillations and the oscillations due to the sources $\Sep_0$ and $\Sep_1$. On the other hand, since the two sources are concentrating close to two different points in $\R^3$, one can guess that they do not interact when $\ep\to 0$. These are the phenomena that the present paper aims at studying quantitatively. We refer to Section 3 for the precise assumptions we need on the sources.\\
\indent In some sense, the sign of the term $-i\aep\ep\uep$ prescribes a radiation condition at infinity for $\uep$. One of the key difficulty in our problem is to follow this condition in the limiting process $\ep\to 0$.\\

We study the high frequency limit in terms of Wigner measures (or semi-classical measures). This is a mean to describe the propagation of quadratic quantities, like the local energy density $|\uep(x)|^2$, as $\ep\to 0$. The Wigner measure $\mu(x,\xi)$ is the energy carried by rays at the point $x$ with frequency $\xi$. These measures were introduced by Wigner~\cite{wig} and then developed by P. G\'erard~\cite{ge} and P.-L. Lions and T. Paul~\cite{lp} (see also the surveys~\cite{bu} and \cite{gmmp}). They are relevant when a typical length $\ep$ is prescribed. They have already proven to be an efficient tool in the study of high frequencies, see for instance \cite{bckp}, \cite{cpr} for Helmholtz equations, P. G\'erard, P.A. Markowich, N.J. Mauser, F. Poupaud~\cite{gmmp} for periodic media, G. Papanicolaou, L. Ryzhik~\cite{pr} for a formal analysis of general wave equations, L. Erd\"os, H.T. Yau~\cite{ey} for an approach linked to statistical physics, and L. Miller~\cite{mil2} for a study in the case with sharp interface.\\

Such problems of high frequency limit of Helmholtz equations have been studied in Benamou, Castella, Katsaounis, Perthame~\cite{bckp} and Castella, Perthame, Runborg~\cite{cpr}. In~\cite{bckp}, the authors considered the case of one point source and a general index of refraction whereas in~\cite{cpr}, they treated the case of a source concentrating close to a general manifold with a constant refraction index. In the present paper, we borrow the methods used in both articles.\\
\indent In the case of one point source, for instance $\Sep_0$ only, with a constant index of refraction, it is proved in \cite{bckp} that the corresponding Wigner measure $\mu_0$ is the solution to the Liouville equation
$$ 0^+\mu_0(x,\xi)+\xi\cdot\nabla_x\mu_0(x,\xi)=
Q_0(x,\xi)=\frac{1}{(4\pi)^2}\delta(x)\delta(|\xi|^2-1)|\widehat{S_0}(\xi)|^2,$$ 
the term $0^+$ meaning that $\mu$ is the outgoing solution given by 
$$ \mu_0(x,\xi)=\int_{-\infty}^0Q_0(x+t\xi,\xi)dt.$$
In particular, the energy source created by $\Sep_0$ is supported at $x=0$. Similarly, the energy source created by the source $\Sep_1$ is supported at $x=q_1$. Thinking of the orthogonality property on Wigner measures, one can guess that the energy source generated by the sum $\Sep_0+\Sep_1$ is the sum of the two energy sources created asymptotically by $\Sep_0$ and $\Sep_1$.\\ 
\indent Indeed, we prove in this paper that the Wigner measure $\mu$ associated with the sequence $(\uep)$ satisfies 
\begin{equation}\label{eqmu}
 0^+\mu(x,\xi)+\xi\cdot\nabla_x\mu(x,\xi)=Q_0+Q_1,
\end{equation}
where $Q_0$ and $Q_1$ are the source terms obtained in \cite{bckp} in the case of one point source. However, our proof does not rest on the mere orthogonality property.\\

Let us now give some details about our proof. Our strategy is borrowed from \cite{bckp}. First, we prove uniform estimates on the sequence of solutions $(\uep)$. We also study the limiting behaviour fo the rescaled solutions $\ep^{\frac{d-1}{2}}\uep(\ep x)$ and $\ep^{\frac{d-1}{2}}\uep(q_1+\ep x)$. The obtention of these first two results is the key difficulty in our paper. It relies on the study of the sequence $(\qep)$ such that
$$-i\aep\ep\qep+\Delta\qep+\qep= S_1\Big(x-\frac{q_1}{\ep}\Big).$$
Using the explicit formula for the Fourier tranform of $\qep$, we prove that $\qep$ is uniformly bounded in a suitable space and that $\qep\to 0$ as $\ep\to 0$ weakly. We would like to point out that our analysis, based on a study in Fourier space, strongly rests on the assumption of a constant index of refraction.\\
Second, our results on the Wigner measure then follow from the properties proved in \cite{bckp}. They are essentially consequences of the uniform bounds on $(\uep)$: we write the equation satisfied by the Wigner transform associated with $(\uep)$, and pass to the limit $\ep\to 0$ in the various terms that appear in this equation. The only difficult (and new) term to handle is the source term.\\
Third, we prove an improved version of the radiation condition of \cite{bckp}. Our argument relies on the observation that $\mu$ is localized on the energy set $\{|\xi|^2=1\}$, a property that was not exploited in \cite{bckp}.\\

The paper is organized as follows. In Section 2, we recall some definitions and state our assumptions. Section 3 is devoted to the proof of uniform bounds on the sequence of solutions $(\uep)$ and of the convergence of the rescaled solutions. Then, in Section 4, we establish the transport equation satisfied by the Wigner measure $\mu$ together with the radiation condition at infinity. In the appendix, we recall the proof of some results established in \cite{bckp} that we use in our paper.

\section{Notations and assumptions}
In this section, we recall the definitions of Wigner transforms and of the $B$, $B^*$ norms introduced by Agmon and H\"ormander~\cite{ah} for the study of Helmholtz equations. Then, we give our assumptions.
 
\subsection{Wigner transform and Wigner measures}
\noindent We use the following definition for the Fourier transform:
$$ \hat{u}(\xi)=(\F_{x\to\xi}u)(\xi)=\frac{1}{(2\pi)^3}\int_{\R^d} e^{-ix\cdot\xi}u(x)dx.$$
For $u,v\in\S(\R^3)$ and $\ep>0$, we define the Wigner transform
\begin{eqnarray*}
\Wep(u,v)(x,\xi)&=&(\F_{y\to\xi})(u\big(x+\frac{\ep}{2}y\big)
\bar{v}\big(x-\frac{\ep}{2}y\big)),\\
\Wep(u)&=&\Wep(u,u).
\end{eqnarray*}
In the sequel, we denote $\Wep=\Wep(\uep)$.\\

If $(\uep)$ is a bounded sequence in $L^2_{loc}(\R^d)$, it turns out that (see~\cite{ge},~\cite{lp}), up to extracting a subsequence, the sequence $(\Wep(\uep))$ converges weakly to a positive Radon measure $\mu$ on the phase space $T^*\R^3=\R^3_x\times\R^3_{\xi}$ called Wigner measure (or semiclassical measure) associated with $(\uep)$:
\begin{equation}
\forall \varphi\in\Co_c^{\infty}(\R^{6}),\ \lim_{\ep\to 0}\langle\Wep(\uep),\varphi\rangle=\int \varphi(x,\xi)d\mu
\end{equation}
 
We recall that these measures can be obtained using pseudodifferential operators. The Weyl semiclassical operator $a^W(x,\ep D_x)$ (or $Op^W_{\ep}(a)$) is the continuous operator from $\S(\R^d)$ to $\S'(\R^d)$ associated with the symbol $a\in\S'(T^*\R^d)$ by Weyl quantization rule
\begin{equation}
 (a^W(x,\ep D_x)u)(x)=\frac{1}{(2\pi)^d}\int_{\R^d_{\xi}}\int_{\R^d_y}a
\left(\frac{x+y}{2},\ep\xi\right)
f(y)e^{i(x-y)\cdot\xi}d\xi dy.
\end{equation}

\noindent We have the following formula: for $u,v\in\S'(\R^d)$ and $a\in\S(\R^d\times\R^d)$,
\begin{equation}\label{wignerpseudo}
\langle\Wep(u,v),a\rangle_{\S',\S}=\langle\bar{v},a^W(x,\ep D_x)\bar{u}\rangle_{\S',\S},
\end{equation}
where the duality brackets $\langle.,.\rangle$ are semi-linear with respect to the second argument. This formula is also valid for $u,v$ lying in other spaces as we will see in Section 3.

\subsection{Besov-like norms}

In order to get uniform (in $\ep$) bounds on the sequence $(\uep)$, we shall use the following Besov-like norms, introduced by Agmon and H\"ormander~\cite{ah}: for $u,f\in L^2_{loc}(\R^3)$, we denote
\begin{eqnarray*}
\|u\|_{B^*}&=&\sup_{j\geq -1}\left(2^{-j}\int_{C(j)}|u|^2dx\right)^{1/2},\\
\|f\|_{B}&=&\sum_{j\geq -1}\Big(2^{j+1}\int_{C(j)}|f|^2dx\Big)^{1/2},
\end{eqnarray*}
where $C(j)$ denotes the ring $\{x\in\R^3/2^j\leq |x|< 2^{j+1}\}$ for $j\geq 0$ and $C(-1)$ is the unit ball.\\

These norms are adapted to the study of Helmholtz operators. Indeed, if $v$ is the solution to 
$$-i\alpha v+\Delta v+v=f $$
where $\alpha>0$, then Agmon and H\"ormander~\cite{ah} proved that there exists a constant $C$ independent of $\alpha$ such that
$$ \|v\|_{B^*}\leq C\|f\|_{B}.$$
Perthame and Vega~\cite{pv} generalised this result to Helmholtz equations with general indices of refraction.\\

\noindent We denote for $x\in \R^3$, $|x|=\sqrt{\sum_{j=1}^{3} x_j^2}$ and $\langle x\rangle = (1+|x|^2)^{1/2}$.\\
 For all $\delta>\frac{1}{2}$, we have
\begin{equation}\label{propb}
\|u\|_{L^2_{-\delta}}:=\|\langle x\rangle^{-\delta}u\|_{L^2}\leq C(\delta)\|u\|_{B^*}.
\end{equation}

We end this section by stating two properties of these spaces that will be useful for our purpose (the reader can find the proofs in \cite{ah}). The first proposition states that, in some sense, we can define the trace of a function in $B$ on a linear manifold of codimension 1.

\begin{prop}\label{pb2}
There exists a constant $C$ such that for all $f\in B$, we have 
$$\int_{\R}\|f(x_1,.)\|_{L^2(\R^{2})}dx_1 \leq C\|f\|_{B}.$$
\end{prop}

The second property gives the stability of the space $B$ by change of variables in Fourier space. 
\begin{prop}\label{pb}
Let $\Omega_1,\ \Omega_2$ be two open sets in $\R^3$, $\psi:\Omega_1\to\Omega_2$ a $\mathcal{C}^2$ diffeomorphism, $\chi\in \mathcal{C}^1_c(\R^3)$. For all $u\in B$, we denote
$$Tu=\mathcal{F}^{-1}\big(\chi(\hat{u}\circ\psi)\big).$$
Then
$$\|Tu\|_{B}\leq C\|\chi\|_{\mathcal{C}^1_b}\|\psi\|_{\mathcal{C}^2_b}\|u\|_{B}.$$
\end{prop}

\subsection{Assumptions}

We are now ready to state our assumptions. Our first assumption, borrowed from \cite{bckp}, concerns the regularizing parameter $\aep>0$.\\
\indent (H1) $ \aep\geq\ep^{\gamma}$ for some $\gamma>0$.\\
This assumption is technical and is used to get a radiation condition at infinity in the limit $\ep\to 0$. Next, in order to get uniform bounds on $\uep$, we assume that the source terms $S_0$ and $S_1$ belong to the natural Besov space that is needed to actually solve the Helmholtz equation (\ref{h}). \\
\indent (H2) $\|S_0\|_B,\| S_1\|_B<\infty$.\\
It turns out that, in order to compute the limit of the energy source, we shall need the stronger assumption\\
\indent (H3) $\langle x\rangle^N S_0\in L^2(\R^3)$ and $\langle x\rangle^N S_1\in L^2(\R^3)$ for some $N>\frac{1}{2}+\frac{3\gamma}{\gamma+1}$.\\

\section{Bounds on solutions to Helmholtz equations}

In this section, we first establish uniform bounds on the sequence $(\uep)$ that will imply estimates on the sequence of Wigner transforms $(\Wep)$. It turns out that we shall also need to compute the limit of the rescaled solutions $\wep_0$ and $\wep_1$ defined below in order to obtain the energy source in the equation satisfied by the Wigner measure $\mu$.\\
Before stating our two results, let us define these rescaled solutions. Following \cite{bckp} and \cite{cpr}, we denote
\begin{equation}\label{defw}
\left\{\begin{array}{lll}
      \wep_0(x)&=&\ep^{\frac{d-1}{2}}\uep(\ep x), \\
      \wep_1(x)&=&\ep^{\frac{d-1}{2}}\uep(q_1+\ep x).
      \end{array}\right.
\end{equation}
They respectively satisfy 
$$\left\{\begin{array}{ccc}
 -i\aep\ep\wep_0+\Delta\wep_0+\wep_0&=&S_0(x)+S_1\big(x-\frac{q_1}{\ep}\big),\\
 -i\aep\ep\wep_1+\Delta\wep_1+\wep_1&=&S_0\big(x+\frac{q_1}{\ep}\big)+S_1(x).
 \end{array}\right.$$
We are ready to state our results on $\uep$, $\wep_0$ and $\wep_1$.\\

\begin{prop}\label{boundu}
Assume $S_0$, $S_1\in B$. Then, the solution $\uep$ to the Helmholtz equation (\ref{h}) satisfies the following bound
$$\|\uep\|_{B^*}\leq C(\|S_0\|_B+\|S_1\|_B),$$
where $C$ is a constant independent of $\ep$.
\end{prop}

\begin{prop}\label{w}
Let $\wep_0$ and $\wep_1$ be the rescaled solutions defined by (\ref{defw}). Then, the sequences $(\wep_0)$ and $(\wep_1)$ are uniformly bounded in $B^*$ and they converge weakly-$\ast$ in $B^*$ to the outgoing solutions $w_0$ and $w_1$ to the following Helmholtz equations
$$
\left\{\begin{array}{l}
   \Delta w_0+w_0=S_0\\
    \Delta w_1+w_1=S_1,
\end{array}\right.
$$
i.e. $w_0$ and $w_1$ are given in Fourier space by 
$$ \widehat{w_j}(\xi)=\frac{-\widehat{S_j}(\xi)}{|\xi|^2-1+i0}
=-\Big(p.v.\big(\frac{1}{|\xi|^2-1}\big)
+i\pi\delta(|\xi|^2-1)\Big)\widehat{S_j}(\xi),\quad j=0,1.$$  
\end{prop}

\rem The Helmholtz equation $\Delta w+w=S $ does not uniquely specify the solution $w$. An extra condition is necessary, for instance the Sommerfeld radiation condition. When the refraction index is constant equal to 1, this condition writes
\begin{equation}\label{somm}
\lim_{r\to\infty}
\frac{1}{r}\int_{S_r}\Big|\frac{\partial{w}}{\partial{r}}+iw\Big|^2d\sigma
=0.
\end{equation}
Such a solution is called an outgoing solution.\\
Alternatively, still assuming that the refraction index is constant, the outgoing solution to the Helmholtz equation may be defined as the weak limit $w$ of the sequence $(w^{\delta})$ such that
$$ -i\delta w^{\delta}+\Delta w^{\delta}+ w^{\delta}=S(x).$$
We point out that the two points of views are equivalent in the case of a constant index of refraction (which is not true for a general index of refraction).\\

We prove the two propositions in the following two sections. As we will see in the proofs, our main difficulties are linked to the rays that are emitted by the source at 0 towards the point $q_1$ (and conversely). Hopefully, the interaction between those rays is "destructive" and not constructive.

\subsection{Proof of Proposition \ref{boundu}}\label{d1}

In the sequel, $C$ will denote any constant independent of $\ep$.\\
The scaling invariance
$$
\|\uep\|_{B^*}\leq \|\wep_0\|_{B^*},
$$
makes it sufficient to prove bounds on $\wep_0$.
\medskip
Since $\wep_0$ is a solution to
$$ -i\aep\ep\wep_0+\Delta\wep_0+\wep_0=S_0(x)+S_1\big(x-\frac{q_1}{\ep}\big)$$
we may decompose $\wep_0=\widetilde{\wep_0}+\qep$, where $\widetilde{\wep_0}$ and $\qep$
satisfy
$$\left\{\begin{array}{ccl}
-i\aep\ep\widetilde{\wep_0}+\Delta\widetilde{\wep_0}+\widetilde{\wep_0}&=&S_0(x),\\
 -i\aep\ep\qep+\Delta\qep+\qep&=&S_1\big(x-\frac{q_1}{\ep}\big).
\end{array}\right.$$
First, we note that the bound $\|\widetilde{\wep_0}\|_{B^*}\leq
C\|S_0\|_B$ is established in Agmon-H\"ormander~\cite{ah} (see also Perthame-Vega~\cite{pv}). Hence, the proof of Proposition \ref{boundu} reduces to the proof of the following lemma.\\

\begin{lemm}\label{q1}
If $\qep$ is the solution to
$$-i\aep\ep\qep+\Delta\qep+\qep=S_1\big(x-\frac{q_1}{\ep}\big)$$
then $\qep$ is uniformly (in $\ep$) bounded in $B^{\star}$:
$$ \|\qep\|_{B^*}\leq C\|S_1\|_B$$
\end{lemm}

\begin{proof} 
We want to prove that
$$\forall v\in B,\quad |\langle\qep,v\rangle|\leq C\|S_1\|_B\|v\|_B.$$
Using Parseval's equality, we write
\begin{equation}\label{sou}
\langle\qep,v\rangle=
\int_{\R^3}\frac{e^{-i\frac{q_1\cdot\xi}{\ep}}\widehat{S_1}(\xi)
\bar{\hat{v}}(\xi)}{-|\xi|^2+1-i\ep\aep}d\xi.
\end{equation}
To estimate this integral, we shall distinguish the values of $\xi$ close to or far from two critical sets: the sphere $\{|\xi|^2=1\}$ (the set where the denominator in (\ref{sou}) vanishes when $\ep\to 0$) and the line
$\{\xi\ \textrm{collinear to}\ x_0\}$ (the set where we cannot apply directly the stationary phase theorem to (\ref{sou})). \\

More precisely, we first take a small parameter $\delta\in]0,1[$, and we distinguish in the integral (\ref{sou}), the contributions due to the values of $\xi$ such that $|\xi^2-1|\geq\delta$ or $|\xi^2-1|\leq\delta$. Let $\chi\in\Co^{\infty}_c(\R)$ be a truncation function such that $\chi(\lambda)=0$ for $|\lambda|\geq 1$. We denote $\chi_{\delta}(\xi)=\chi\big(\frac{|\xi|^2-1}{\delta}\big)$. We accordingly decompose
\begin{eqnarray*}
\langle\qep,v\rangle&=&
\int_{\R^3}\frac{e^{-i\frac{q_1}{\ep}\cdot\xi}\widehat{S_1}(\xi)
\bar{\hat{v}}(\xi)\chi_{\delta}(\xi)}{-|\xi|^2+1-i\ep\aep}d\xi
+\int_{\R^3}\frac{e^{-i\frac{q_1}{\ep}\cdot\xi}\widehat{S_1}(\xi)
\bar{\hat{v}}(\xi)(1-\chi_{\delta}(\xi))}{-|\xi|^2+1-i\ep\aep}d\xi\\
&=& I^{\ep}+II^{\ep}. 
\end{eqnarray*}

First, since the denominator is not singular on the support of $\chi_{\delta}$, we easily bound the first part with the $L^2$ norms 
$$
|I^{\ep}|
\leq\frac{\|\chi\|_{L^{\infty}}}{\delta}\|\widehat{S_1}\|_{L^2}\|\hat{v}\|_{L^2},
$$
and using $B\hookrightarrow L^2$, we obtain the desired bound
\begin{equation}
|I^{\ep}| \leq C \|S_1\|_{B}\|v\|_{B}.
\end{equation}

Let us now study the second part $II^{\ep}$ where the denominator is singular.
Up to a rotation, we may assume $q_1=|q_1|e_1$, where $e_1$ is the first vector of the canonical base. We make the polar change of variables 
$$\xi=\left\{\begin{array}{l}
              r\sin\theta\cos\varphi\\
              r\sin\theta\sin\varphi\\
              r\cos\theta
             \end{array}\right..
$$

\rem In order to make the calculations easier, we write this paper in dimension equal to 3, but the proof would be similar in any dimension $d\geq 3$.\\

\noindent Hence, 
$q_1\cdot\xi=|q_1|r\sin\theta\cos\varphi$,
and we get
$$II^{\ep}=\int\frac{e^{-i\frac{|q_1|}{\ep}r\sin\theta\cos\varphi}}{-r^2+1-i\ep\aep}
\big(\widehat{S_1}\bar{\hat{v}}(1-\chi_{\delta})\big)(\xi(r,\theta,\varphi))r^2\sin\theta drd\theta d\varphi$$
Now, we distinguish the contributions to the integral $d\theta d\phi$ linked to the values close to, or far from, the critical direction $\{\theta=\frac{\pi}{2},\ \varphi=0\}$ (which corresponds to the case $\{\xi\ \textrm{collinear to}\ q_1\}$). To that purpose, let $\eta>0$ be a small parameter and denote 
$$K=\left\{(r,\theta, \phi)\ \bigg| \ 1-\chi\Big(\frac{r^2-1}{\delta}\Big)\neq 0, \chi\left(\frac{\theta-\frac{\pi}{2}}{\eta}\right)\neq 0, \chi\Big(\frac{\varphi}{\eta}\Big)\neq 0\right\}.$$ 
Then $K$ is a compact set. Let $k\in\Co^{\infty}_c$ be such that $(1-\chi_{\delta})k(\theta,\varphi)$ is a localization function on $K$. We write
\begin{eqnarray*}
II^{\ep}&=&\int\frac{e^{-i\frac{|q_1|}{\ep}r\sin\theta\cos\varphi}}
{-r^2+1-i\ep\aep}
(\widehat{S_1}\bar{\hat{v}})(\xi(r,\theta,\varphi))
(1-\chi_{\delta}(r))k(\theta,\varphi)r^2\sin\theta drd\theta d\varphi\\
& +&\int\frac{e^{-i\frac{|q_1|}{\ep}r\sin\theta\cos\varphi}}
{-r^2+1-i\ep\aep}
(\widehat{S_1}\bar{\hat{v}})(\xi(r,\theta,\varphi))
(1-\chi_{\delta}(r))(1-k(\theta,\varphi))r^2\sin\theta drd\theta d\varphi\\
II^{\ep}&=& III^{\ep}+IV^{\ep}.
\end{eqnarray*}

To estimate the contribution $III^{\ep}$, we apply the stationary phase method. We denote $\alpha=\theta-\frac{\pi}{2}$. The phase function is
$g_r(\alpha,\varphi)=r\cos\alpha\cos\varphi$ so
$$ \begin{array}{rcccll}
  \displaystyle\frac{\partial g_r}{\partial \alpha}&=&-r\sin\alpha\cos\varphi&=&0&\ \textrm{at} \ (\alpha, \varphi)=(0,0),\\
  \displaystyle\frac{\partial g_r}{\partial \varphi}&=&-r\cos\alpha\sin\varphi&=&0&\ \textrm{at}\  (\alpha, \varphi)=(0,0),
   \end{array}$$
and the Hessian at the point $(0,0)$ is
$$D^2g_r(0,0)=\left(\begin{array}{cc}
                     -r&0\\
                     0&-r
                      \end{array}\right),
$$
which is invertible at any point in $K$. Since $K$ is a compact set, we can apply the Morse lemma: there exists a finite covering $(\Omega_j)_{j=1,n}$ ($n\in\N$) of $K$ such that on each set $\Omega_j$, there exists a $\Co^{\infty}$ change of variables 
$(\alpha,\varphi)\mapsto(\alpha_j,\varphi_j)$ such that
$$
g_r(\alpha,\varphi)=r-r\frac{\alpha_j^2}{2}-r\frac{\varphi_j^2}{2}.
$$
Moreover, we can write $(1-\chi_{\delta})k=\sum_{j=1}^n\chi_j$ where $\chi_j\in\Co^{\infty}_c$ and $supp(\chi_j)\subset\Omega_j$. Then, we make the
changes of variables $\alpha_j'=\sqrt{\frac{r}{2}}\alpha_j$,
$\varphi_j'=\sqrt{\frac{r}{2}}\varphi_j$. Finally, we decompose $\chi_j=\chi_j^1\chi_j^2$. Thus, we obtain, for the
contribution $III^{\ep}$, the formula
\begin{equation}\label{terme3}
III^{\ep}=\sum_{j=1}^n
\int\frac{e^{i\frac{x_1}{\ep}(-r+\alpha_j'^2+\varphi_j'^2)}}{-r+1+i\ep\aep}
\widehat{T_j^1S_1}(r,\alpha_j',\varphi_j')
\overline{\widehat{T_j^2v}}(r,\alpha_j',\varphi_j')drd\alpha_j'
d\varphi_j',
\end{equation}
where
\begin{eqnarray*}
T_j^1S_1&:=&\mathcal{F}\Big((\chi_j^1\widehat{S_1})\circ
\xi(r,\alpha(\alpha_j,\varphi_j),\varphi(\alpha_j,\varphi_j))\Big),\\
T_j^2v&:=&\mathcal{F}\bigg(\frac{-r+1+i\ep\aep}{-r^2+1-i\ep\aep}(\chi_j^2\hat{v})
\circ\xi(r,\alpha(\alpha_j,\varphi_j),\varphi(\alpha_j,\varphi_j))\\
&&\qquad\qquad\times\frac{2}{r}\Big|\frac{\textrm{d}\xi}
{\textrm{d}(r,\alpha,\varphi)}
\Big|\Big|\frac{\textrm{d}(\alpha,\varphi)}{\textrm{d}(\alpha_j,\varphi_j)}\Big|
\bigg).
\end{eqnarray*}
As a first step, using Proposition \ref{pb}, we directly get $T_j^1S_1\in B$ with
$$\|T_j^1S_1\|_B\leq C\|S_1\|_B.$$
As a second step, we study $T_j^2v$. Since for $r$ close to $1$,
$$\Big|\frac{-r+1+i\ep\aep}{-r^2+1-i\ep\aep}\Big|\leq 1,$$
we recover, from Proposition \ref{pb}, 
$$T_j^2v\in B\qquad\textrm{and}\qquad\|T_j^2v\|_B\leq C\|v\|_B.$$
Now, we apply Parseval's equality with respect to the r variable in the formula (\ref{terme3}) 
\begin{eqnarray*}
III^{\ep}&=&\sum_{j=1}^n
\int\frac{e^{i\frac{|q_1|}{\ep}(\alpha_j'^2+\varphi_j'^2)}}{-r+1+i\ep\aep}
\widehat{T_j^1S_1(.-\frac{|q_1|}{\ep},.,.)}\overline{\widehat{T_j^2v}}drd\alpha_j' d\varphi_j'\\
&=&\int
e^{i\frac{|q_1|}{\ep}(\alpha_j'^2+\varphi_j'^2)}\1_{\{\rho>0\}}e^{-(\ep\aep-i)t}
\mathcal{F}_{r\to\rho}(\widehat{T_j^1S_1(.-\frac{|q_1|}{\ep},.,.)})
(\rho-t,\alpha_j',\varphi_j')\\
&&\qquad\qquad\qquad\times\mathcal{F}_{r\to\rho}(\overline{\widehat{T_j^2v}})
(\rho,\alpha_j',\varphi_j')dt
d\rho d\alpha_j' d\varphi_j'.
\end{eqnarray*}
where $\1_{\{\rho>0\}}$ denotes the characteristic function of the set $\{\rho>0\}$.\\
Hence, we obtain 
\begin{eqnarray*}
|III^{\ep}|&\leq &\sum_{j=1}^n\left(\int\|\mathcal{F}_{r\to\rho}
(\widehat{T_j^1S_1(.-\frac{|q_1|}{\ep},.,.)})
(\rho))\|_{L^2}d\rho\right)\\
&&\qquad\qquad\times\left(\int\|\mathcal{F}_{r\to\rho}
(\widehat{T_j^2v})(\rho)\|_{L^2}d\rho\right)\\
|III^{\ep}|&\leq & \sum_{j=1}^n\Big(\int\|T_j^1S_1(\rho-\frac{|q_1|}{\ep})\|_{L^2}d\rho\Big)
\Big(\int\|T_j^2v(\rho)\|_{L^2}d\rho\Big)\\
|III^{\ep}|&\leq &\sum_{j=1}^n \Big(\int\|T_j^1S_1(\rho)\|_{L^2}d\rho\Big)
\Big(\int\|T_j^2v(\rho)\|_{L^2}d\rho\Big)\\
|III^{\ep}|&\leq & C\sum_{j=1}^n\|T_j^1S_1\|_B\|T_j^2v\|_B\\
|III^{\ep}|&\leq & C\|S_1\|_B\|v\|_B,
\end{eqnarray*}
which is the desired estimate.\\

We are left with the part $IV^{\ep}$, which corresponds to the directions $\xi$ that are not collinear to $q_1$. We denote $K'$ the support of $(1-\chi_{\delta})(1-k)$ which is a compact set. In $K'$, we can choose as new independent variables
$$ \eta_1=-q_1\cdot\xi,\quad \eta_2=|\xi|^2-1.$$
More precisely, since
$$\frac{\textrm{d}(\eta_1,\eta_2)}{\textrm{d}\xi}=\left(\begin{array}{c}
                                                             -q_1\\
                                                             2\xi
                                                             \end{array}\right)$$
is of maximal rank 2, there exists a finite covering $(\Omega_j')_{j=1,m}$ ($m\in\N$) of $K'$ such that in $\Omega_j'$, we can make the
change of variables $\xi\mapsto\eta$. As before, we denote
$\chi_j'=\chi_j^3\chi_j^4$ some localization functions on $\Omega_j'$ such that $ (1-\chi_{\delta})(1-k)=\sum_{j=1}^m\chi_j'$. Thus, for $j=1,\dots,m$,
$$
\int\frac{e^{-i\frac{q_1}{\ep}\cdot\xi}}{-|\xi|^2+1+i\ep\aep}
\widehat{S_1}\overline{\hat{v}}\chi_j' d\xi
=\int\frac{e^{i\frac{\eta_2}{\ep}}}{-\eta_1+i\ep\aep}
(\widehat{S_1}\overline{\hat{v}}\chi_j')(\xi(\eta))
\Big|\frac{\textrm{d}\xi}{\textrm{d}\eta}\Big|d\eta.
$$
If we denote
\begin{eqnarray*}
T_j^3S_1&:=&\mathcal{F}^{-1}((\chi_j^3\widehat{S_1})\circ\xi),\\
T_j^4v&:=&\mathcal{F}^{-1}\big((\chi_j^4\hat{v})\circ\xi)
\big|\frac{\textrm{d}\xi}{\textrm{d}\eta}\big|\big),
\end{eqnarray*}
and if $\mathcal{F}_1$ denotes the Fourier transform with respect to the $\eta_1$ variable, Parseval's equality with respect to $\eta_1$ gives
\begin{eqnarray*}
\bigg|\int\frac{e^{-i\frac{q_1}{\ep}\cdot\xi}}{-|\xi|^2+1+i\ep\aep}
\widehat{S_1}\overline{\hat{v}}\chi_j'
d\xi\bigg|
&=&(2\pi)^d\bigg|\int\chi_{\{t>0\}}e^{-\ep\aep t}(\mathcal{F}_1^{-1}(\widehat{T_j^3S_1}))(x_1-t)\\
&&\qquad\times(\mathcal{F}_1^{-1}
(\overline{\widehat{T_j^4v}}))(x_1)e^{i\eta_2 /\ep}dtdx_1d\eta_2 d\eta_3 \bigg| \\
&\leq & C\|S_1\|_B\|v\|_B.
\end{eqnarray*}
Summing over $j$, we obtain
$$ |IV^{\ep}|\leq C \|S_1\|_B\|v\|_B,$$
which ends the proof of the bound 
$$ |\langle\qep,v\rangle|\leq C \|S_1\|_B\|v\|_B.$$
\end{proof}

\subsection{Proof of Proposition \ref{w}}

We prove the result for the sequence $(\wep_0)$, the convergence of the sequence $(\wep_1)$ can be obtained similarly. As we did in the proof of Proposition \ref{boundu}, we write $\wep_0=\widetilde{\wep_0}+\qep$. Since $\widetilde{\wep_0}$ is the solution to a Helmholtz equation with constant index of refraction and fixed source, it converges weakly-$\ast$ to the outgoing solution $w_0$ to $ \Delta w +w=S_0$. Hence, it suffices to show the following result.\\

\begin{lemm}\label{q2}
If $\qep\in B^{\star}$ is the solution to
$$-i\aep\ep\qep+\Delta\qep+\qep=S_1\big(x-\frac{q_1}{\ep}\big)$$
then $\qep\to 0$ in $B^{\star}$.
\end{lemm}
\begin{proof}
The proof of this result requires two steps (using a density argument):
\begin{enumerate}
\item for $v\in B$, we have the bound $\big|\langle\qep,v\rangle\big|\leq C\|S_1\|_{B}\|v\|_{B}$
\item if $S_1$ and $v$ are smooth, then $\langle\qep,v\rangle\to 0$.
\end{enumerate}
The first point is exactly the result in Lemma \ref{q1}. It remains to prove the convergence in the smooth case (the second point above).\\
We write
$$
\langle\qep,v\rangle=\int_{\R^3}\frac{e^{-i\frac{q_1}{\ep}\cdot\xi}\widehat{S_1}(\xi)
\bar{\hat{v}}(\xi)}{-|\xi|^2+1-i\ep\aep}d\xi.
$$

We are thus left with the study of 
\begin{equation}
R_{\ep}(\psi)=\int_{\R^3}\frac{e^{-i\frac{q_1}{\ep}\cdot\xi}
\psi(\xi)}{-|\xi|^2+1-i\ep\aep} d\xi
\end{equation}
where
$\psi=\widehat{S_1}\overline{\widehat{v}}$ belongs to $\mathcal{S}(\R^3)$.\\

As in the proof of Lemma \ref{q1}, we distinguish the
contributions of various values of $\xi$. We shall use exactly the same partition, according to the values of $\xi$ close to, or far from, the sphere $|\xi|=1$ and collinear or not to $q_1$. We shall use the same notations for the various truncation functions. \\
We first separate the contributions of $\xi$ such that $|\xi^2-1|\leq \delta$ and $|\xi^2-1|\geq \delta$
\begin{eqnarray*}
R_{\ep}(\psi)&=&
\int_{\R^3}\frac{e^{-i\frac{q_1}{\ep}\cdot\xi}
\psi(\xi)\chi_{\delta}(\xi)}{-|\xi|^2+1-i\ep\aep}d\xi
+\int_{\R^3}\frac{e^{-i\frac{q_1}{\ep}\cdot\xi}\psi(\xi)
(1-\chi_{\delta}(\xi))}{-|\xi|^2+1-i\ep\aep}d\xi\\
&=& I^{\ep}+II^{\ep}. 
\end{eqnarray*}

In the support of $\chi_{\delta}$, since the denominator is not singular, we can apply the non stationary phase method.\\
Since $q_1\neq 0$, we may assume $q_1^1\neq 0$ and we have
\begin{eqnarray*}
I^{\ep}
&=&\frac{\ep}{iq_1^1}\int_{\R^3}e^{-i\frac{q_1}{\ep}\cdot\xi}
\partial_{\xi_1}\bigg(
\frac{\psi(\xi)\chi_{\delta}(\xi)}{-|\xi|^2+1-i\ep\aep}\bigg)d\xi\\
{}&=&\frac{\ep}{iq_1^1}\int_{\R^3}e^{-i\frac{q_1}{\ep}\cdot\xi}\bigg(
\frac{\partial_{\xi_1}(\psi(\xi)\chi_{\delta}(\xi))}{-|\xi|^2+1-i\ep\aep}
-\frac{2\psi(\xi)\chi_{\delta}(\xi)\xi_1}
{(-|\xi|^2+1-i\ep\aep)^2}\bigg).
\end{eqnarray*}
Hence, we obtain the bound
\begin{eqnarray*}
|I^{\ep}|&\leq &\frac{\ep}{|q_1^1|} \int_{\R^3}\Big(\frac{1}{\delta}|\partial_{\xi_1}(\chi\psi)|
+\frac{2}{\delta^2}|\xi_1\chi\psi|\Big)d\xi.
\end{eqnarray*}
Since $\partial_{\xi_1}(\chi\psi)$ and $\xi_1\chi\psi$ belongs to $\mathcal{S}$, we have, as $\ep\to 0$,
$$I^{\ep}\to 0.$$

Let us now study the second term $II^{\ep}$. We use the same changes of variables as in Section \ref{d1}. It leads to the following formula
$$
II^{\ep}=\sum_{j=1}^n
 \int\frac{e^{-i\frac{q_1}{\ep}(r-\alpha_j'^2-\varphi_j'^2)}}{-r+1+i\ep\aep}
\widetilde{\chi_j}(r,\alpha_j',\varphi_j')
\widetilde{\psi}(r,\alpha_j',\varphi_j')drd\alpha_j'd\varphi_j',$$
where
\begin{eqnarray*}
\widetilde{\chi_j}(r,\alpha_j',\varphi_j')&=&
\chi_j\circ\xi(r,\alpha(\alpha_j',\varphi_j'),
\varphi(\alpha_j',\varphi_j'))\frac{2(-r+1+i\ep\aep)}{r(-r^2+1-i\ep\aep)}
\bigg|\frac{\textrm{d}(\alpha,\varphi)}{\textrm{d}(\alpha_j,\varphi_j)}\bigg|,\\
\widetilde{\psi}(r,\alpha_j',\varphi_j')&=&
\psi\circ\xi(r,\alpha(\alpha_j',\varphi_j'),
\varphi(\alpha_j',\varphi_j')),
\end{eqnarray*}
are still smooth functions that are bounded independently from $\ep$.\\
Using Parseval's inequality with respect to the variables $(\alpha_j',\varphi_j')$ for each integral, we obtain the bound
$$|II^{\ep}|\leq 
C\ep\sum_{j=1}^n 
\left|\int\frac{e^{-i\frac{|q_1|}{\ep}r}e^{-i\ep(\lambda_j^2+\mu_j^2)}}
{-r+1+i\ep\aep}\mathcal{F}_{\lambda_j,\mu_j}
(\widetilde{\chi_j}\widetilde{\psi})drd\lambda_j
d\mu_j\right|.$$

To obtain the convergence of $II^{\ep}$, it remains to study an integral of the following type
$$
\int_{|r-1|\leq\delta}\frac{e^{-i\frac{|q_1|}{\ep}r}w(r)}{-r+1+i\ep\aep}dr,\
\textrm{where}\ w\in\mathcal{S}.
$$
This is done in the following lemma.
\begin{lemm}\label{l}
$\forall w\in \mathcal{S},\ \forall \theta\in\,(0,1)$, we have
$$\int_{|r|\leq\delta}\frac{e^{-i\frac{|q_1|}{\ep}r}w(r)}{-r+i\ep\aep}dr
=-i\pi w(0)+O_{\ep\to 0}(\ep^{-\theta}).
$$
\end{lemm}

\noindent Using this lemma, we readily get the estimate
\begin{equation}
|III^{\ep}|
\leq C\ep^{1-\theta}\quad \forall \theta \in\,(0,1),
\end{equation}
which proves that $III^{\ep}\to 0 $ as $\ep\to 0$.\\

\noindent There remains to give the\\
\noindent{\em Proof of Lemma \ref{l}.} We write
\begin{eqnarray*}
\int_{-\delta}^{\delta}\frac{e^{-i\frac{|q_1|}{\ep}r}w(r)}{-r+i\ep\aep}dr
&=&\int_{-\delta}^{\delta}\frac{e^{-i\frac{|q_1|}{\ep}r}w(r)}{r^2+(\ep\aep)^2}(r-i\ep\aep)dr\\
&=&-i\ep\aep\int_{-\delta}^{\delta}\frac{e^{-i\frac{|q_1|}{\ep}r}w(r)}
{r^2+(\ep\aep)^2}dr
+\int_{-\delta}^{\delta}e^{-i\frac{|q_1|}{\ep}r}
\frac{rw(r)}{r^2+(\ep\aep)^2}dr\\
&=&I+II.
\end{eqnarray*}
We have
$$
I=-i\int_{-\frac{\delta}{\ep\aep}}^{\frac{\delta}{\ep\aep}}\frac{e^{-i|q_1|\aep
y}w(\ep\aep y)}{y^2+1}dy \to -i\pi w(0),$$ 
and
$$ II=\int_{-\delta}^{\delta}\big(e^{-i\frac{|q_1|}{\ep}r}w(r)-w(0)\big)
\frac{r}{r^2+(\ep\aep)^2}dr+
\int_{-\delta}^{\delta}w(0)\frac{r}{r^2+(\ep\aep)^2}dr\ 
.
$$
The last term vanishes because the integrand is odd. Moreover, using the smoothness of $w$, we easily obtain that for all $\theta\in\,(0,1)$,
$$\big|e^{-i\frac{|q_1|}{\ep}r}w(r)-w(0)\big|\leq
C_{\theta}\Big(\frac{r}{\ep}\Big)^{\theta}$$
Thus,
$$\bigg|\int_{-\delta}^{\delta}\big(e^{-i\frac{|q_1|}{\ep}r}w(r)-w(0)\big)
\frac{r}{r^2+(\ep\aep)^2}dr\bigg|
\leq
\frac{C}{\ep^{\theta}}\int_{-\delta}^{\delta}|r|^{\theta-1}dr$$ and
the result is proved. \fin\\

We are left with the study of $IV^{\ep}$.
We use the same change of variables as in Section \ref{d1}. 
\begin{eqnarray*}
IV^{\ep}&=&\sum_{j=1}^m
\int\frac{e^{-i\frac{q_1}{\ep}\cdot\xi}}{-|\xi|^2+1+i\ep\aep}
\psi(\xi)\chi_j(\xi)d\xi\\
&=&\sum_{j=1}^m\int\frac{e^{i\frac{\eta_2}{\ep}}}{-\eta_1+i\ep\aep}
(\psi\chi_j)(\xi(\eta))
\Big|\frac{\textrm{d}\xi}{\textrm{d}\eta}\Big|d\eta \\
&=&i\ep\sum_{j=1}^m
\int\frac{e^{i\frac{\eta_2}{\ep}}}{-\eta_1+i\ep\aep}
\partial_{\eta_1}\Big((\psi\chi_j)(\xi(\eta))
\Big|\frac{\textrm{d}\xi}{\textrm{d}\eta}\Big|\Big)d\eta.
\end{eqnarray*}
The integral obviously converges with respect to all the variables except $\eta_1$. It remains to prove the convergence with respect to the $\eta_1$ variable, i.e. the convergence of
$$\int\frac{\phi(\eta)}{-\eta_1+i\ep\aep}d\eta_1,$$
where 
$$\phi=\partial_{\eta_1}\Big((\psi\chi_j)(\xi(\eta))
\Big|\frac{\textrm{d}\xi}{\textrm{d}\eta}\Big|\Big)$$
is smooth and compactly supported with respect to $\eta$. It is a consequence of the fact that the distribution $(x+i0)^{-1}$ is well-defined on $\R$ by
$$\frac{1}{x+i0}=v.p.(\frac{1}{x})-i\pi\delta(x).$$
We conclude that $IV^{\ep}\to 0$ and $\langle \qep,v\rangle\to 0$ as $\ep\to 0$. 
\end{proof}

\section{Transport equation and radiation condition on $\mu$}

In this section, we state and prove our results on the Wigner measure associated with $(\uep)$. Since we established the uniform bounds on $(\uep)$ and the convergence of $(\wep_0)$, $(\wep_1)$, these results now essentially follows from the results proved in \cite{bckp}. We first prove bounds on the sequence of Wigner transforms $(\Wep)$ that allow us to define a Wigner measure $\mu$ associated to $(\uep)$. Then, we get the transport equation satisfied by $\mu$ together with the radiation condition at infinity, which uniquely determines $\mu$.

\subsection{Results}

\begin{theo}\label{tt1}
Let $S_0,\ S_1\in B$ and $\lambda>0$. The sequence $(\Wep)$ is
bounded in the Banach space $X_{\lambda}^{\star}$ and up to extracting a subsequence, it converges weak-$\star$ to a positive and locally bounded measure $\mu$ such that
\begin{equation}\label{bornew}
\sup_{R>0}\frac{1}{R} \int_{|x|<R}\int_{\xi\in
\R^3}\mu(x,\xi)\,dxd\xi \leq C(\|S_0\|_B+\|S_1\|_B)^2.
\end{equation}
The Banach space $X^*_{\lambda}$ is defined as the dual space of the set $X_{\lambda}$ of functions $\hat{\varphi}(x,\xi)$ such that $\varphi(x,y):=\F_{\xi\to y}(\hat{\varphi}(x,\xi))$ satisfies
\begin{equation}\label{xlambda}
\int_{\R^d}\sup_{x\in\R^d}(1+|x|+|y|)^{1+\lambda}|\varphi(x,y)|dy <\infty.
\end{equation}
\end{theo}

\begin{theo}\label{tt2}
Assume (H1), (H2), (H3). Then the Wigner measure
$\mu$ asssociated with $(\uep)$ satisfies the following transport equation
\begin{equation}\label{eqf2}
\xi\cdot\nabla_{x}\mu=\frac{1}{(4\pi)^2}\Big(\delta(x)|\widehat{S_0}(\xi)|^2
+\delta(x-q_1)|\widehat{S_1}(\xi)|^2\Big)\delta(|\xi|^2-1):=Q(x).
\end{equation}
Moreover, $\mu$ is the outgoing solution
to the equation (\ref{eqf2}) in the following sense: for all test function $R \in \Co^{\infty}_c(\R^6)$, if we denote $ g(x,\xi)=\int_{0}^{\infty}R(x-\xi t,\xi)dt$, then
\begin{equation}\label{rad}
\int_{\R^6}R(x,\xi)d\mu(x,\xi)=-\int_{\R^6}Q(x,\xi)g(x,\xi)dxd\xi.
\end{equation}
\end{theo}

\rem Here the support of the test function $R$ contains 0, contrary to \cite{bckp}.

\subsection{Proof of Theorem \ref{tt1}}
This theorem, that is proved in \cite{bckp}, is a consequence of the uniform estimate on the sequence $(\uep)$ in the space $B^*$ obtained in Proposition \ref{boundu}. We observe that for any $\lambda>0$,
\begin{equation}
\|\langle x\rangle^{-\frac{1}{2}-\lambda} \uep(x)\|_{L^2}\leq C\|\uep\|_{B^*}\leq C\|f\|_{B},
\end{equation}
hence, for any function $\varphi$ satisfying (\ref{xlambda}), we have 
\begin{eqnarray*}
&&|\langle \Wep(\uep),\hat{\varphi}\rangle|\\
&&\leq  \int_{\R^{6}}\frac{|\uep|(x+\frac{\ep}{2}y)|\overline{\uep}|(x-\frac{\ep}{2}y)}{\langle x+\frac{\ep}{2}y \rangle ^{\frac{1}{2}+0}
\langle x-\frac{\ep}{2}y \rangle ^{\frac{1}{2}+0}}\langle x+\frac{\ep}{2}y \rangle ^{\frac{1}{2}+0}\langle x-\frac{\ep}{2}y\rangle ^{\frac{1}{2}+0}|\varphi|(x,y)dxdy\\
&&\leq  C\|f\|_{B}^{2}\int_{\R^3} \sup_{x\in \R^3} \langle |x|+|y|\rangle ^{1+0}|\varphi(x,y)|dy.
\end{eqnarray*}

\noindent So $(\Wep(\uep))$ is bounded in $X^*_{\lambda}$, $\lambda >0$. We deduce that, up to extracting a subsequence, $(\Wep(\uep))$ converges weak-$\ast$ to a nonnegative measure $\mu$ satisfying  
\begin{equation}\label{est2}
|\langle \mu,\hat{\varphi}\rangle|
 \leq C\|f\|_{B}^{2} \int_{\R^3} \sup_{x\in \R^3} \langle|x|+|y|\rangle ^{1+0}|\varphi(x,y)|dy.
\end{equation}
We refer for instance to Lions, Paul~\cite{lp} for the proof of the nonnegativity of $\mu$.\\
The bound (\ref{bornew}) is obtained using the following family of functions
$$
\varphi^R_{\mu}(x,y)=\frac{1}{\mu^{3/2}}e^{-|y|^2/\mu}\frac{1}{R}\chi(\langle
x\rangle\leq R)
$$
and letting $\mu\to 0$, $R\to\infty$.\hfill\fin\\

\subsection{Proof of the transport equation \ref{eqf2}}

This section is devoted to the proof of the transport equation satisfied by $\mu$. We first write the transport equation satisfied by $\Wep$ in a dual form. Then, we study the convergence of the source term (the convergence of the other terms is obvious). Finally, choosing an appropriate test function in the limiting process, we get the radiation condition at infinity satisfied by $\mu$. Proving first a localization property, we improve the radiation condition proved in \cite{bckp}.

\subsubsection{Transport equation satisfied by $\Wep$}

$\Wep$ satisfies the following equation
\begin{equation}\label{eqW}
\aep\Wep+\xi\cdot\nabla_x\Wep=
\frac{i\ep}{2}\mathcal{I}m\;\Wep(\fep,\uep):=Q^{\ep}.
\end{equation}
This equation can be obtained writing first the equation satisfied by $$\vep(x,y)=\uep\big(x+\frac{\ep}{2}y\big)
\overline{\uep}\big(x-\frac{\ep}{2}y\big).$$ 
From the equality
$$ \nabla_{y}\cdot\nabla_{x}\vep=
\frac{\ep}{2}\big[\Delta\uep(x+\frac{\ep}{2}y)\overline{\uep}(x-\frac{\ep}{2}y)
-\Delta\overline{\uep}(x-\frac{\ep}{2}y)\uep(x+\frac{\ep}{2}y)\big],
$$
we deduce 
$$\aep\vep+i\nabla_{y}\cdot\nabla_{x}\vep+\frac{i}{2\ep}
\big[n^2(x+\frac{\ep}{2}y)-n^2(x-\frac{\ep}{2}y)\big]\vep=\sigma_{\ep}(x,y),$$
where
$$\sigma_{\ep}(x,y):=
\frac{i\ep}{2}\big[\Sep(x+\frac{\ep}{2}y)\overline{\uep}(x-\frac{\ep}{2}y)
-\overline{\Sep}(x-\frac{\ep}{2}y)\uep(x+\frac{\ep}{2}y)\big].
$$
After a Fourier transform, we obtain the equation (\ref{eqW}).\\
Then we write the dual form of this equation. Let $\psi\in\S(\R^6)$, we have
\begin{equation}\label{duale}
\aep\langle\Wep,\psi\rangle-\langle\Wep,\xi\cdot\nabla_x\psi\rangle
=\langle\Qep,\psi\rangle.
\end{equation}
By the definition of the Wigner measure $\mu$, we get
$$ \aep\langle\Wep,\psi\rangle\to 0\qquad
\textrm{and}\qquad 
\langle\Wep,\xi\cdot\nabla_x\psi\rangle\to\langle\mu,\xi\cdot\nabla_x\psi\rangle.$$
Hence we are left with the study of the source term $\langle\Qep,\psi\rangle$.

\subsubsection{Convergence of the source term}
In order to compute the limit of the source term in (\ref{eqW}), we develop
$$\langle\Qep,\psi\rangle = \frac{i\ep}{2}\mathcal{I}m\Big( \langle\Wep(\Sep_0,\uep),\psi\rangle+\langle\Wep(\Sep_1,\uep),\psi\rangle\Big).$$
Thus, the result is contained in the following proposition. \\

\begin{prop}\label{termeQ}
The sequences $\big(\ep\Wep(\Sep_0,\uep)\big)$ and  $\big(\ep\Wep(\Sep_1,\uep)\big)$ are bounded in $\S'(\R^6)$ and for all $\psi\in\S(\R^6)$, we have
\begin{eqnarray}\label{source}
\lim_{\ep\to 0}\ep\langle\Wep(\Sep_0,\uep),\psi\rangle_{\S',\S}& =& \frac{1}{(2\pi)^3}\int_{\R^3} \overline{\widehat{w_0}}(\xi)\widehat{S_0}(\xi) \psi(0,\xi)d\xi,\\
\lim_{\ep\to 0}\ep\langle\Wep(\Sep_1,\uep),\psi\rangle_{\S',\S}& =& \frac{1}{(2\pi)^3}\int_{\R^3} \overline{\widehat{w_1}}(\xi)\widehat{S_1}(\xi) \psi(q_1,\xi)d\xi,
\end{eqnarray}
where $w_0$  and $w_1$ are defined in Proposition \ref{w}.
\end{prop}

\noindent Using Proposition \ref{termeQ}, we readily get 
\begin{eqnarray*}
\lim_{\ep\to 0}\langle\Qep,\psi\rangle &= &
\frac{i}{2(2\pi)^3}\mathcal{I}m\left( \int_{\R^3} \overline{\widehat{w_0}}(\xi)\widehat{S_0}(\xi) \psi(0,\xi)d\xi+\int_{\R^3} \overline{\widehat{w_1}}(\xi)\widehat{S_1}(\xi) \psi(q_1,\xi)d\xi\right)\\
&=& \frac{1}{(4\pi)^2}\bigg( \int_{\R^3} |\widehat{S_0}(\xi)|^2\delta(\xi^2-1) \psi(0,\xi)d\xi\\
&&\qquad\qquad\qquad\qquad +\int_{\R^3} |\widehat{S_1}(\xi)|^2\delta(\xi^2-1) \psi(q_1,\xi)d\xi\bigg),
\end{eqnarray*}
which is the result in Theorem \ref{tt2}. \fin\\

\noindent Let us now prove Proposition \ref{termeQ}.\\
{\em Proof of Proposition \ref{termeQ}.} The two terms to study being of the same type, we only consider the first one in our proof. Let $\psi\in\S(T^*\R^d)$ and $\varphi(x,y)=\F^{-1}_{y\to\xi}(\psi(x,\xi))$, then we have
\begin{eqnarray*}
\ep\langle\Wep(\Sep_0,\uep),\psi\rangle_{\S',\S}
&=&\ep\int\Sep_0\big(x+\frac{\ep}{2}y\big)\overline{\uep}\big(x-\frac{\ep}{2}y\big)
\varphi(x,y)dxdy \\
&=&\int S_0(x)\overline{\wep_0}(x+y)\varphi(\ep(x+\frac{y}{2}),y)dxdy.
\end{eqnarray*}
Hence, using that $\psi\in\S(\R^{2d})$, we get
\begin{eqnarray*}
\big|\ep\langle\Wep(\Sep_0,\uep),\psi\rangle_{\S',\S}\big|
&\leq &  C\int \langle x\rangle^N |S_0(x)|\frac{|\wep_0(x+y)|}{\langle x+y \rangle^{\beta}}\frac{\langle x+y\rangle^{\beta}}{\langle x\rangle^N\langle y\rangle^k}dxdy\\
&\leq & C\|\langle x\rangle^N S_0\|_{L^2}\|\wep_0\|_{B^*}
\int_{\R^3_y}\sup_{x\in\R^3}\frac{\langle x+y\rangle^{\beta}}{\langle x\rangle^N\langle y\rangle^k}dy
\end{eqnarray*}
for any $k\geq 0$ and $\beta >1/2$, upon using the Cauchy-Schwarz inequality in $x$.\\
Then, we distinguish the cases $|x|\leq |y|$ and $|x|\geq |y|$ : the term stemming from the first case gives a contribution which is bounded by $C\int\frac{dy}{\langle y\rangle^{k-\beta}}$ and the second contribution is bounded by 
$C\int\frac{dy}{\langle y\rangle^{k}}$. So, upon choosing $k$ large enough, we obtain that
 $$|\ep\langle\Wep(\Sep_0,\uep),\psi\rangle_{\S',\S}|\leq C\|\langle x\rangle^N S_0\|_{L^2}\|\wep\|_{B^*}.$$
Now, in order to compute the limit (\ref{source}), we write
\begin{eqnarray*}
\ep\langle\Wep(\Sep_0,\uep),\psi\rangle
&=&\int S_0(x)\overline{\wep_0}(x+y)\bigg(\varphi\Big(\ep\Big(x+\frac{y}{2}\Big),y\Big)
-\varphi(0,y)\bigg)dxdy\\
&&+\int \overline{\wep_0}(x)S_0(x-y)\varphi(0,y)dxdy\\
&=&I_{\ep}+II_{\ep}.
\end{eqnarray*}
Reasonning as above, we readily get that $\lim_{\ep\to 0}I_{\ep}=0$. For the second term, we have
$$II_{\ep}=\int \overline{\wep_0}(x)\big(S_0\ast\varphi(0,.)\big)(x)dx,$$
hence, since $\wep_0$ converges weakly-$\ast$ in $B^*$, it suffices to prove that $S_0\ast\varphi(0,.)$ belongs to $B$. We denote $\phi=\varphi(0,.)$. We have, for $\beta>1/2$,
\begin{eqnarray*}
\|S_0\ast\phi\|_B &\leq & C \|S_0\ast\phi\|_{L^2_\beta}
=C\int \langle x\rangle^{\beta}|S_0\ast\phi(x)|^2dx \\
&\leq & C\|\phi\|_{L^1}\int \langle x\rangle^{\beta}|S_0|^2\ast|\phi|(x)dx
\end{eqnarray*}
where we used the Cauchy-Schwarz inequality. Hence, we get
$$
\|S_0\ast\phi\|_B \leq 
C \|\langle x\rangle^N S_0\|_{L^2}\int_{\R^3_y}\sup_{x\in\R^3} \frac{\langle x+y\rangle^{\beta}}{\langle x\rangle^N\langle y\rangle^k}dy,
$$
for any $k$. As before, this integral converges. Thus, we have established that $S_0\ast\widehat{\psi}(0,.)$ belongs to $B$, which implies that 
$$II_{\ep}\to\int S_0(x)\overline{w_0}(x+y)\widehat{\psi}(0,y)dxdy.$$
\fin

\subsection{Proof of the radiation condition (\ref{rad})}

It remains to prove that $\mu$ satisfies the weak radiation condition (\ref{rad}). 

\subsubsection{Support of $\mu$}

In order to prove the radiation condition without restriction on the test function $R$ (as assumed in \cite{bckp}), we first prove a localization property on the Wigner measure $\mu$. This property is well-known when $\uep$ satisfies a Helmholtz equation without source term. It is still valid here thanks to the scaling of $\Sep$.
\begin{prop} Under the hypotheses (H1), (H2), (H3), the Wigner measure $\mu$ satisfies
$$
{\rm supp}(\mu)\subset\{(x,\xi)\in \R^6/\ |\xi|^2=1\}.
$$
\end{prop}

\begin{proof}
Let $\phi\in\Co^{\infty}_c(\R^6)$ and $\ph=\phi^W(x,\ep D_x)$. 
Let us denote $\hep=-\ep^2\Delta-1$. Since $\uep$ satisfies the Helmholtz equation (\ref{h}), we have
\begin{equation}\label{h2}
 i\aep\ep\uep+\hep\uep=\ep^2\Sep.
\end{equation}
Moreover, $\hep$ is a pseudodifferential operator with symbol $|\xi|^2-1 $. By pseudodifferential calculus, $ \ph\hep=Op^W_{\ep}(\phi(x,\xi)(|\xi|^2-1))+O(\ep)$ so, using the definition of the measure $\mu$, we get that
\begin{eqnarray*}
\lim_{\ep\to 0}\,(\ph\hep\uep,\uep)&=&
\lim_{\ep\to 0}\,(Op^W_{\ep}(\phi(x,\xi)(|\xi|^2-1))\uep,\uep)\\
&=& \int \phi(x,\xi)(|\xi|^2-1))d\mu.
\end{eqnarray*}
Using the equation (\ref{h2}), we write
$$
(\ph\hep\uep,\uep) =\ep^2(\ph\Sep,\uep)-i\aep\ep(\ph\uep,\uep)
=\ep^2(\Wep(\Sep,\uep),\phi)-i\aep\ep(\ph\uep,\uep). 
$$
On the first hand, Proposition \ref{termeQ} gives that $\lim_{\ep\to 0}\ep^2(\Wep(\Sep,\uep),\phi)= 0$. On the other hand, $(\ph\uep,\uep)$ is bounded so $\lim_{\ep\to 0}\aep\ep(\ph\uep,\uep)= 0$. Therefore, for any $\phi\in\Co^{\infty}_c(\R^6)$, we have $\int \phi(|\xi|^2-n^2(x))d\mu=0$, so ${\rm supp}(\mu)\subset\{|\xi|^2=1\}$. 
\end{proof}

\subsubsection{ Proof of the condition (\ref{rad})}

Using the previous localization property, in order to prove the radiation condition (\ref{rad}), one may only use test functions $R\in\Co^{\infty}_c(\R^6)$ such that ${\rm supp}(R)\subset\R^6\backslash \{\xi=0\}$.\\
Let $R$ be such a test function. We associate with $R$ the solution $\gep$ to
$$-\aep\gep+\xi\cdot\nabla_x\gep=R(x,\xi).$$
By duality, we have
$$\langle\Qep,\gep\rangle=\langle\Wep,R\rangle,$$
so that it sufffices to establish the following two convergences:
\begin{equation}\label{l1}
\lim_{\ep\to 0}\langle\Qep,\gep\rangle=\langle Q,g\rangle,
\end{equation}
\begin{equation}\label{l2}
\lim_{\ep\to 0}\langle\Wep,R\rangle=\langle f,R\rangle,
\end{equation}
where $Q$ et $g$ are defined in Theorem \ref{tt2}.\\
As before, since $R\in X_{\lambda}$ for any $\lambda>0$, the limit (\ref{l2}) follows from the weak-$\ast$ convergence of $\Wep$ in $X_{\lambda}^{\star}$.\\
On the other hand,
\begin{eqnarray}\label{cv}
\langle\Qep,\gep\rangle&=&\mathcal{I}m\int_{\R^6}\overline{S_0}(x)\wep_0(x+y)
\widehat{\gep}(\ep[x+\frac{y}{2}],y)dxdy\\\nonumber
&&\qquad\qquad+\mathcal{I}m\int_{\R^6}\overline{S_1}(x)\wep_1(x+y)\widehat{\gep}
(q_1+\ep[x+\frac{y}{2}],y)dxdy,
\end{eqnarray}
so $\langle\Qep,\gep\rangle$ is the sum of two terms of the same type. Such a term has been studied in \cite{bckp}, where the following result is proved.

\begin{prop}\label{condrad} 
Assume $(\wep)$ is bounded in $B^*$ and that $(\wep)$ converges weakly-$\ast$ in $B^*$ to $w_0$. Assume $S_0$ satisfy (H3). Let $R\in\Co^{\infty}_c(\R^6)$ be such that $supp(R)\subset\R^6\backslash \{\xi=0\}$. Let $\gep$ be the solution to 
$$ -\aep\gep+\xi\cdot\nabla_x\gep=R(x,\xi)$$
and $g(x,\xi)=\int_0^{\infty}R(x+t\xi,\xi)dt$. Then, we have
$$\lim_{\ep\to 0}\int_{\R^6}\overline{S_0}(x)\wep_0(x+y)
\widehat{\gep}(\ep[x+\frac{y}{2}],y)dxdy=
\frac{1}{(2\pi)^3}\int_{\R^3} \overline{\widehat{S_0}}(\xi)\widehat{w_0}(\xi)g(0,\xi)d\xi.$$
\end{prop}

\begin{proof} 
The proof of this result is written in the appendix.
\end{proof}

\noindent Using the proposition above together with Proposition \ref{w}, we get that  
\begin{eqnarray*}
&&\lim_{\ep\to 0}
\langle\Qep,\gep\rangle\\
&&=\mathcal{I}m\bigg(\frac{1}{(2\pi)^3}\int_{\R^3} \overline{\widehat{S_0}}(\xi)\widehat{w_0}(\xi)g(0,\xi)d\xi
+\frac{1}{(2\pi)^3}\int_{\R^3} \overline{\widehat{S_1}}(\xi)\widehat{w_1}(\xi)g(q_1,\xi)d\xi\bigg)\\
&&= \frac{1}{(4\pi)^2}\left(\int_{\R^3} |\widehat{S_0}(\xi)|^2·\delta(\xi^2-1)g(0,\xi)d\xi
+\int_{\R^3} |\widehat{S_1}(\xi)|^2\delta(\xi^2-1)g(q_1,\xi)d\xi\right).
\end{eqnarray*}
Thus, the radiation condition (\ref{rad}) is proved.

\appendix
\section{Proof of Proposition \ref{condrad}}
In the sequel, we denote 
$$G^{\ep}=\lim_{\ep\to 0}\int_{\R^6}\overline{S_0}(x)\wep_0(x+y)
\widehat{\gep}(\ep[x+\frac{y}{2}],y)dxdy.$$

\subsection{Bounds on $G^{\ep}$}\label{borne}
 
In order to study $G^{\ep}$, we need a preliminary result on the test function $\gep$.
\begin{lemm}\label{lemg}
Let $R\in\Co^{\infty}_c(\R^6\setminus\{\xi=0\})$. We denote $\gep$ the solution to
\begin{equation}
-\aep \gep+\xi.\nabla_{x}\gep=R(x,\xi).
\end{equation}
It is given by the explicit formula
\begin{equation}
 \gep(x,\xi)=-\int_{0}^{\infty} \exp(-\aep|\xi|^{-1}s)\frac{1}{|\xi|}R(x-\frac{\xi}{|\xi|}s,\xi)ds.
\end{equation}
Then we have the estimate
\begin{equation}
\forall M\geq 0, \quad |\widehat{\gep}(x,y)|\leq C\frac{\langle x \rangle ^{M}\wedge\aep^{-M}}{\langle y \rangle ^{M}},
\end{equation}
\noindent where $.\wedge.$ denotes the infimum of two numbers, and $C$ is a constant depending on $M$ and $R$.
\end{lemm}

\begin{proof} Let $a$ be a multiindex such that $|a|\leq M$. We denote
$\omega=\frac{\xi}{|\xi|}$. We write
\begin{eqnarray}\label{gep}
y^{a}\widehat{\gep}(x,y)\!\!&=&\!\!\mathcal{F}_{\xi\to y}\Big(\int_{0}^{\infty}(i\partial_{\xi})^{a}\bigg[e^{-\aep|\xi|^{-1}s}\frac{1}{|\xi|}R(x-\omega s,\xi)\bigg]ds\Big) \nonumber\\
&=&\!\!\!\!\displaystyle\int_{\R^3}d\xi e^{-i\xi.y}\!\int_{s=0}^{+\infty}\!\!ds\sum_{b,c,d,e,f}C_{(a,b,c,d,e,f)}
e^{-\aep|\xi|^{-1}s}(i\partial_{\xi})^{b}(-\aep|\xi|^{-1}s)\nonumber\\
&&
\!\!\times(-\aep|\xi|^{-1}s)^{f}(i\partial_{\xi})^{c}(-s\omega)(i\partial_{x})^{d}
(i\partial_{\xi})^{e}\bigg(\frac{1}{|\xi|}R\bigg)(x-\omega s,\xi).
\end{eqnarray}
Using that\\
$\bullet$ $R\in \Co_{0}^{\infty}(\R^6\backslash\{|\xi|=0\})$ so\\
\indent $\bullet$ there exists $r_0,\ A,\ B>0$ such that $supp(R)\subset \{|x|\leq r_0\}\times\{A\leq|\xi|\leq B\},$\\
\indent $\bullet$ $R$ and its derivatives belong to $L^1(\R^6)$.\\
\noindent $\bullet$ $ |(i\partial_{\xi})^{b}(-\aep|\xi|^{-1}s)|\leq Cs, \quad \forall\ |\xi|\geq A$.\\
$\bullet$ $|(i\partial_{\xi})^{c}(-s\omega)|\leq Cs, \quad \forall\ |\xi|\geq A$.\\
$\bullet$ in the integral above, $s\in
[|x|-r_0,|x|+r_0]$, so for $|x|$ large enough, we can use the equivalence $s\sim|x|\sim\langle x\rangle$ (where we dente for $a,\ b>0$, $a\sim b$ if $\exists\ c_1,\ c_2>0/\
c_1a<b<c_2a$), we get in (\ref{gep}),
$$
|y^{a}\widehat{\gep}(x,y)|\leq C\langle x \rangle ^{M}\exp(-\frac{\aep}{B}\langle x \rangle).
$$
The desired estimate follows.
\end{proof}
%\hfill\fin\\

\noindent Using this lemma, we estimate
\begin{eqnarray*}
|G^{\ep}|&\leq& \Big|\int_{\R^6}\wep_0(x+y)S_0(x)\widehat{\gep}(\ep(x+\frac{y}{2},y)dxdy\Big|\\
&\leq&C\int_{\R^6}\frac{|\wep_0(x+y)|}{\langle x+y\rangle^{\frac{1}{2}+0}}\langle x+y\rangle^{\frac{1}{2}+0}\langle x\rangle^{N_1}|S_0(x)|\langle x\rangle^{-N_1}\frac{\langle\ep(x+\frac{y}{2})\rangle^{M}\wedge\aep^{-M}}{\langle y \rangle^{M}}\\
&\leq&C\|\wep_0\|_{B^{\star}}\|\langle x\rangle^{N_1}S_0(x)\|_{L^2}\\
&&\qquad\qquad\times\int_{\R^3}\sup_{x\in\R^3}\langle
x+y\rangle^{\frac{1}{2}+0}\langle
x\rangle^{-N_1}\frac{\langle\ep(|x|+|y|)\rangle^{M}\wedge\aep^{-M}}{\langle
y \rangle^{M}}dy.
\end{eqnarray*}
Now, we prove that the integral
$$\mathcal{I}^{\ep}=\int_{\R^3}\sup_{x\in\R^3}\langle
x+y\rangle^{\frac{1}{2}+0}\langle
x\rangle^{-N_1}\frac{\langle\ep(|x|+|y|)\rangle^{M}\wedge\aep^{-M}}{\langle
y \rangle^{M}}dy$$
is bounded uniformly with respect to $\ep$.\\
Let us define the following three subsets in $\R^6$
\begin{equation}\begin{array}{l}
A_{\ep}=\{(x,y)\in \R^6/|x|\geq |y|\}, \quad B_{\ep}=\{|x|\leq|y|,|\ep^{1-0} y|\leq 1\}, \\
C_{\ep}=\{|x|\leq| y|,|\ep^{1-0} y|\geq 1\},
\end{array}
\end{equation}
where $\ep^{1-0}$ means $\ep^{1-\delta}$ with $\delta >0$ sufficiently small.\\
\underline{If $(x,y)\in A_{\ep}$}, then 
$$\langle
x+y\rangle^{\frac{1}{2}+0}\langle
x\rangle^{-N_1}\frac{\langle\ep(|x|+|y|)\rangle^{M}\wedge\aep^{-M}}{\langle
y \rangle^{M}}
\leq C\langle
x\rangle^{-N_1+\frac{1}{2}+0}\Big(\langle\ep x\rangle^{M}
\wedge\aep^{-M}\Big)\langle y \rangle^{-M}.$$
Now, we distinguish the relative size of $\langle \ep
x\rangle$ and $\ep^{\gamma}$:\\
\indent $\bullet$ if $\langle\ep x\rangle\geq\ep^{-\gamma}$, we have
$$
\langle x\rangle^{\frac{1}{2}+0-N_1}\Big(\langle\ep x\rangle^{M}\wedge\aep^{-M}\Big)\leq\langle x\rangle^{\frac{1}{2}+0-N_1}\ep^{-\gamma M}\leq\ep^{-\gamma(\frac{1}{2}+0-N_1+M)}.
$$
\indent $\bullet$ if $\langle\ep x\rangle\leq\ep^{-\gamma}$, we have
$\langle x\rangle\leq\frac{1}{\ep}\langle\ep x\rangle$, hence we get
$$
\langle x\rangle^{\frac{1}{2}+0-N_1}\Big(\langle\ep x\rangle^{M}\wedge\aep^{-M}\Big)\leq\langle x\rangle^{\frac{1}{2}+0-N_1}\ep^{-\gamma M}\leq\ep^{-\gamma M}\ep^{-(\gamma+1)(\frac{1}{2}+0-N_1)}.
$$
Now, we choose $1/2+0$ such that $N_1<\frac{1}{2}+0$. Then, we get the following bound for the contribution of the set $A_{\ep}$ to $\mathcal{I}_{\ep}$
$$C\ep^{-\gamma M-(\gamma+1)(1/2+0-N_1)}\int_{\R^3}\langle y\rangle^{-M}dy.$$
Since
$N_1>\frac{3\gamma}{\gamma+1}+\frac{1}{2}$, we can choose $M>3$
such that this contribution is uniformly bounded with respect to $\ep$.\\
\underline{If $(x,y)\in B_{\ep}$}, then $|\ep y|\leq 1$ so we obtain 
$$\langle
x+y\rangle^{\frac{1}{2}+0}\langle
x\rangle^{-N_1}\frac{\langle\ep(|x|+|y|)\rangle^{M}\wedge\aep^{-M}}{\langle
y \rangle^{M}}
\leq C\langle
y\rangle^{-M+\frac{1}{2}+0}.$$
Thus, the corresponding contribution to $\mathcal{I}^{\ep}$ is bounded by $ C\int_{\R^3}\langle y\rangle^{-M+\frac{1}{2}+0}dy$ which is convergent.\\
\underline{If $(x,y)\in C_{\ep}$}, then
$$\langle
x+y\rangle^{\frac{1}{2}+0}\langle
x\rangle^{-N_1}\frac{\langle\ep(|x|+|y|)\rangle^{M}\wedge\aep^{-M}}{\langle
y \rangle^{M}}
\leq C
\langle y\rangle ^{\frac{1}{2}+0}\frac{\langle \ep y\rangle ^{M}\wedge\aep^{-M}}{\langle y\rangle ^{M}}.$$
Since $|y|\geq \ep^{-1+0}$, if we denote $z=\ep y$, we have
$|z|\geq\ep^{\delta}$, for some $\delta >0$, which implies that
$$
\ep\langle \frac{z}{\ep}\rangle \geq \frac{\ep^{\delta}}{\sqrt{2}}\langle z\rangle.
$$
Thus, we have
\begin{eqnarray*}
\mathcal{I}_{\ep}&\leq& C \ep^{-3}\int_{\R^3}\langle \frac{z}{\ep}\rangle ^{-M+\frac{1}{2}+0}(\langle z\rangle ^{M}\wedge\ep^{-\gamma M})dz\\
&\leq & C \ep^{M(1-\delta)-4}\int_{\R^3}\langle z\rangle ^{-M+\frac{1}{2}+0}(\langle z\rangle ^{M}\wedge\ep^{-\gamma M})dz
\end{eqnarray*}
where we used the hypothesis (H1) $\aep\geq\ep^{\gamma}$.\\
We distinguish two cases, according to the relative size of $\langle z\rangle$ and $\ep^{-\gamma}$. We have
\begin{eqnarray*}
\int_{\langle z\rangle \geq
\ep^{-\gamma}}\langle z\rangle ^{\frac{1}{2}+0} \frac{\langle z\rangle
^{M}\wedge\ep^{-\gamma M}}{\langle z\rangle ^{M}}dz
&\leq &\ep^{-M\gamma}\int_{\langle z\rangle \geq \ep^{\gamma}}\langle z\rangle ^{\frac{1}{2}+0-M}dz\\
&\leq & C \ep^{-\gamma(\frac{3}{2}+0)},
\end{eqnarray*}
and
\begin{eqnarray*}
\int_{\langle z\rangle \leq
\ep^{\gamma}}\langle z\rangle ^{\frac{1}{2}+0}\frac{\langle z\rangle
^{M}\wedge\ep^{-\gamma M}}{\langle z\rangle ^{M}}dz
&\leq &\int_{\langle z\rangle \leq \ep^{\gamma}}\langle z\rangle ^{\frac{1}{2}+0}dz\\
&\leq& C \ep^{-\gamma(\frac{3}{2}+0)}.
\end{eqnarray*}
Hence, we get
$$ \mathcal{I}_{\ep}\leq C \ep^{M(1-\delta)-4-\gamma(\frac{3}{2}+0)}.$$
To conclude, we choose $M>\frac{4+\frac{3}{2}\gamma}{1-\delta}$, which gives that $\mathcal{I}_{\ep}$ is uniformly bounded with respect to $\ep$.

\subsection{Convergence of $G^{\ep}$}
We decompose $G^{\ep}$ in the following way
\begin{eqnarray*}
G^{\ep}&=&\int_{\R^6}\wep_0(x+y)S_0(x)\big(
\widehat{\gep}(\ep(x+\frac{y}{2}),y)-\widehat{\gep}(0,y)\big)dxdy \\
&&+\int_{\R^6}\wep_0(x+y)S_0(x)
\big(\widehat{\gep}(0,y)-\widehat{g}(0,y)\big)dxdy\\
&&+\int_{\R^6}\wep_0(x+y)S_0(x)\widehat{g}(0,y)dxdy \\
&=&I_{\ep}+II_{\ep}+III_{\ep}
\end{eqnarray*}

Using the same method as in Section \ref{borne}, we prove that $I_{\ep},\ II_{\ep}\to 0$. Then, we may write
$$
III_{\ep}=\int_{\R^6}\wep_0(x)\big(S_0\ast\widehat{g}(0,.)\big)(x)dx.
$$
Moreover, we established in the proof of Proposition \ref{termeQ} that if $\phi$ is rapidly decreasing at infinity, then  $S_0\ast\phi\in B$. Hence, since $(\wep_0)$ converges weakly-$\ast$ in $B^*$ to $w_0$, we get
$$ \lim_{\ep\to 0} III_{\ep}=\int_{\R^6}w_0(x)\big(S_0\ast\widehat{g}(0,.)\big)(x)dx,$$
which ends the proof.\\

{\em Acknowledgement.} I would like to thank my advisor Fran\c cois Castella for having guided this work.

\end{document}